\documentclass[preprint,11pt]{elsarticle}
\usepackage{enumitem}
\usepackage{amssymb}
\usepackage{eurosym}
\usepackage{color}
\usepackage{bm}
\usepackage{booktabs}
\usepackage{enumitem}
\usepackage{multirow}
\usepackage{longtable}
\usepackage{arydshln}
\usepackage[top=2.2cm,bottom=2.5cm,left=2.5cm,right=2.5cm]{geometry}
\usepackage{hyperref}
\usepackage{mathtools}
\usepackage{algorithm,algpseudocode}
\usepackage{subcaption}
\usepackage{tikz}
\usepackage{pgfplots}
\usetikzlibrary{patterns}
\usetikzlibrary{arrows}
\usepackage{amsmath}
\usepackage{etoolbox}
\usepackage{setspace}
\usepackage{etoolbox}

\usepackage{amsthm,amsmath,eurosym}

\newtheorem{definition}{Definition}
\newtheorem{theorem}{Theorem}
\newtheorem{proposition}{Proposition}
\newtheorem{example}{\emph{Example}}

\newtheorem{corollary}{Corollary}
\newcommand{\R}{\mathbb{R}}

\newcommand{\pt}{\text{ }\forall\text{ }}
\newcommand{\tq}{\text{ }:\text{ }}

\newcommand{\Z}{\mathbb{Z}}

\newcommand{\cmi}{\mathcal{M}_n(\mathcal{I}(\R))}

\newcommand{\cn}{\mathcal{N}}
\newcommand{\ci}{\mathcal{I}}

\pgfplotsset{compat=1.18} 

\journal{European Journal of Operational Research}

\usepackage{lineno, blindtext}

\makeatletter
\def\ps@pprintTitle{%
  \let\@oddhead\@empty
  \let\@evenhead\@empty
  \def\@oddfoot{\reset@font\hfil\thepage\hfil}
  \let\@evenfoot\@oddfoot
}
\makeatother

\begin{document}

\begin{frontmatter}
    \title{Interval-valued Value Functions from Uncertain Preferences}

    \author[JAEN1,cor]{Diego {\sc  Garc{\'i}a-Zamora}}\ead{dgzamora@ujaen.es}
    \author[CEGIST]{Jos{\'e} {\sc Rui~Figueira}}\ead{figueira@tecnico.ulisboa.pt }
    
    \address[JAEN1]{Department of Mathematics, University of Ja{\'e}n, 23071 Ja{\'e}n, Spain}
    \address[CEGIST]{CEGIST, Instituto Superior T\'{e}cnico,  Universidade de Lisboa, Portugal}
    \cortext[cor]{Corresponding author at: Department of Mathematics, University of Ja{\'e}n, 23071 Ja{\'e}n, Spain}
    \begin{abstract}
    \noindent The construction of numerical value scales (or priority values) is a recurrent topic in decision-aiding research. However, in real contexts, uncertainty and limited cognitive precision often lead decision-makers to provide interval judgments rather than exact values. In this scenario, even though obtaining a numerical value scale from interval preferences could be feasible, it implies a loss of information and an oversimplification of the input information. This paper proposes a general framework for deriving interval-valued value scales from interval-valued pairwise comparisons. This implies addressing fundamental challenges regarding the elicitation, the determination of representative value functions, and the preservation of properties such as monotonicity and interpretability. We start by presenting a new definition for the consistency of an interval-valued preference relation. We show that when consistency holds, interval-valued value functions can be determined, and we establish conditions for their uniqueness. Furthermore, we show that our proposal extends and unifies classical models such as fuzzy preference relations, Saaty's preference relations, and the Deck of Cards Method. Methodologically, we develop optimization models that provide procedures to guide decision-makers towards consistency when contradictions arise. The framework provides a coherent and interpretable foundation for constructing monotonic interval value functions, bridging classical and interval-valued preference models, and enhancing robustness in decision-making.    
    \end{abstract}
    \vspace{0.25cm}
\end{frontmatter}

\vfill\newpage

\section{Introduction}

\noindent  In the classical literature, the construction of a value scale implies assigning numerical values to a set of objects in such a way that the associated values reflect the degree of difference between the objects \cite{dyer1979}. Although there are plenty of studies devoted to obtaining priority values \cite{Saaty2008}, they neglect that, in real decision-making scenarios, uncertainty and limited cognitive precision often lead decision-makers to feel more comfortable providing interval judgments rather than exact values. However, when the input information is uncertain, it is not sensible to expect that the resulting value scale has a numerical nature, since this implies an oversimplification of the problem and a loss of information \cite{WOS:001016726700001,GARCIAZAMORA2024FSS}.

Therefore, the main goal of this paper is to provide a coherent and explainable framework that enables the derivation of an interval-valued value scale from interval-valued pairwise judgements. This shift from crisp to interval information poses significant challenges: How should the elicitation process be designed? Under what conditions can an interval-valued value function be determined? And how can we ensure that such functions preserve desirable properties such as monotonicity or interpretability? Addressing these questions requires a robust theoretical basis for managing interval information \cite{WOS:000897578400019}.

In preference modeling and decision analysis, consistency becomes an essential concept when the elicitation process is performed through pairwise comparisons. In this regard, consistency (or transitivity) is a necessary condition for the quality and robustness of the information provided by decision-makers, as it guarantees that the given information does not present contradictions. In classical settings, consistency has been defined as an algebraic equation that a preference table must satisfy. This is the case, for instance, in Saaty’s Analytic Hierarchy Process (AHP) \cite{Saaty2008}, Fuzzy Preference Relations (FPRs) \cite{Orlovski1978}, or the Deck of Cards Method \cite{corrente2021}. Furthermore, consistency is also associated with the possibility of determining representative priority values (or a monotonic value scale) from the preferences given by the decision-maker. 

For that reason, in this paper, we propose a general framework for working with pairwise interval comparisons where consistency is rigorously defined for interval information, leading to the construction of interval-valued value functions that are univocally determined \cite{WOS:001536866000001,WOS:001536900700001}. Our definition ensures that if the provided interval judgments satisfy the consistency conditions, then there exists a monotonic interval value function representing them. We also provide unicity conditions.

In addition, we show that our proposal generalizes many classical preference models based on pairwise comparisons. In this sense, it naturally unifies previous definitions of consistency across different paradigms: (i) in the case of FPRs, our framework recovers the usual additive transitivity conditions, (ii) in Saaty’s preference relations (SPRs), it reproduces the multiplicative consistency condition, and (iii) in the Deck of Cards Method, it coincides with the constructed value scale. 

 From a methodological perspective, we further develop optimization models that serve two complementary purposes. First, they allow us to test whether a given set of interval judgments is consistent under our definition. Second, they provide a systematic procedure to guide the decision-maker towards consistent interval information when inconsistencies arise. In this way, the proposed models play both a diagnostic and a prescriptive role, facilitating interactive decision support while respecting the bounded uncertainty inherent in human judgments. Additionally, we analyze the conditions of existence and uniqueness of the solution, giving the exact formula to compute optimal solutions analytically.

The remainder of this paper is as follows. Section \ref{sec:IVPR} introduces a rigorous definition of consistency for interval judgments. Section \ref{sec:scales}, shows that such a definition leads to the determination of monotonic interval value functions. Section \ref{sec:classical} is devoted to proving that our framework extends and unifies classical approaches such as FPRs, SPRs, and the Deck of Cards Method. Section \ref{sec:opt} introduces several optimization models that detect and correct inconsistencies, assisting decision-makers in producing reliable interval information. Section \ref{sec:example}, presents two examples regarding the computational processes to derive value scales by focusing on the Deck of Cards method, thus extending it to handle interval information. Finally, Section \ref{sec:conclusion} concludes the paper.

%

\section{Preliminaries}\label{sec:prelim}
\noindent In this section, we briefly introduce the basic notions that sustain this research. 

\subsection{Interval arithmetic}
\noindent Given an interval $I\subseteq\R$, we denote by $\ci (I)$ the set of intervals contained in $I$. Given an interval $z\in\ci(\R)$, we will write $z=[z^-,z^+]$ to denote its lower and upper bounds, while $\ell(z)=z^+-z^-$ stands for the length of the interval, i.e., the uncertainty associated with $z$. Since intervals allow quantifying the uncertainty on the measure of a certain quantity, researchers developed arithmetic rules to enable computations with interval values \cite{WOS:000535181800002,SANTANA202044}. The arithmetic rules that we consider in this paper are those related to the classical rules for operating with sets \cite{Wang2009}:
\begin{itemize}[label={--}]
    \item \textit{Interval Addition}: \begin{equation*}
        z_1+z_2=\{x_1+x_2\tq x_1\in z_1, x_2\in z_2\}=[z_1^++z_2^-,z_1^++z_2^+], \pt z_1,z_2\in \ci(\R),
    \end{equation*}
    \item \textit{Interval Multiplication by scalar}:
    \begin{equation*}
        \lambda z=\{\lambda x\tq x\in z\}=[\lambda z^-,\lambda z^+], \pt \lambda\in\R, z\in\ci(\R).
    \end{equation*}
    \item \textit{Partial order}:\begin{equation*}
        z_1\leq_0 z_2\iff z_1^-\leq z_2^-\text{ and }z_1^+\leq z_2^+, \pt z_1,z_2\in \ci(\R)
    \end{equation*}
\end{itemize}

Although other authors may use different symbols to denote addition and multiplication by scalars, we keep the same notation as for real numbers since there is no room for a misunderstanding here, and thus we can facilitate the reading of the manuscript. Furthermore, we aim at emphasizing that the interval values can be treated as real numbers, in the sense that they can be ordered and aggregated with the corresponding methodology. Indeed, we may see $\ci(I)$ as a cardinal scale whose elements are uncertain.

\subsection{Deck of Cards Method}
\noindent Within the many elicitation methodologies, pairwise comparisons allow the decision-makers to provide evaluations by focusing on one pair of objects at a time, which may relax the cognitive burden in the process. There are several preference structures based on pairwise comparisons, such as FPRs or SPRs \cite{Orlovski1978, Saaty2008}. 

In this manuscript, we will pay special attention to the pairwise comparison tables obtained through the Deck of Cards method \cite{corrente2021}. Such a methodology is a socio-technical approach that is presented to the decision-maker as a process in which he/she expresses his/her preference intensity regarding two objects to be compared by using a deck of cards. Of note, this approach has recently been applied to elicit fuzzy numbers from the decision-makers \cite{GARCIAZAMORAEJOR}. Below, we describe the basic interaction protocol. Before starting, the decision-maker is provided with a set of cards representing the objects (scale levels, alternatives, criteria, etc.) and a set of blank cards (as many as necessary). The process for gathering the preference information is as follows:

\begin{enumerate}
    \item Rank the objects from the worst to the best, if the same have the same ``attractiveness'' or ``importance'', they should be grouped in the same ranking position.  
    \item The more or less ``attractiveness'' between rank positions is modeled through the insertion of blank cards,
        \begin{enumerate}
            \item Zero blank cards between two rank positions do not mean that they have the same attractiveness but that the difference is minimal, i.e., equal to the unit (this concept will be introduced later on).
            \item One blank card means that the difference in attractiveness between the two objects is twice the unit. 
            \item And so on.
        \end{enumerate}
    \item For computational purposes, in order to build a cardinal scale, the decision-maker should provide two reference levels, and the worst and the best as considered as such reference levels.
\end{enumerate}

The complete process is detailed in \cite{corrente2021}, whereas a summary and application to model fuzzy numbers can be found in \cite{GARCIAZAMORAEJOR}. In addition, this paper also explains the complete process in Section \ref{sec:example}, in which we extend the methodology to be applied when the decision maker provides his/her opinions using an uncertain number of cards. 


\section{Interval-valued Preference Relations}\label{sec:IVPR}
\noindent This section aims to extend the idea of a preference relation to the interval settings. Of course, our goal is to obtain a meaningful extension in which properties such as reciprocity and consistency can be satisfied by generalizing their classical analogs for real numbers.

Although later on we will discuss some transformations to relate our structure with classical preference tables, such as FPRs or SPRs, in this paper, we work in the following set
\begin{equation*}
    \cmi=\left\{Z=\left(z_{ij}\right)_{\substack{1\leq i\leq n\\1\leq j\leq n}}\tq z_{ij}\in\ci\left(\R\right), \, \pt i,j=1,...,n\right\},
\end{equation*}
which stand for the order $n$ matrices whose elements are intervals contained in $\R$. Note that $z_{32}$ is an interval that represents the numbers that we use to model the preference intensity between the objects $3$ and $2$, expressed through an uncertain number of units. In general, the greater $z_{32}$ is with respect to the interval $[0,0]$, the higher the utility of $3$ over the utility of $2$. In contrast, the lower $z_{32}$ is with respect to the interval $[0,0]$, the lower the utility of $3$ over the utility of $2$. If $z_{32}$ is close to $[0,0]$, it means that $3$ and $2$ have a similar utility. Of course, these degenerate intervals are difficult to obtain. What we may have in practice is, for example, that $z_{32}=[-0.05,0.075]$, which can be interpreted as the utility of objects $3$ and $2$ are similar, but we do not know with complete certainty if one of them is better than the other. We remark here that this discussion has involved terms such as ``utility'' or ``similar to'' when talking about intervals. At this stage, we are using these terms as vague notions, not in a strict mathematical sense. We will come back to this issue later.

A key property of classical PRs is that the elements on the main diagonal of the matrix represent the notion of neutrality. Wherever the scale is multiplicative (SPR) or additive (FPR), when an object is compared to itself, the intensity should be a numeric value representing the indifference. In our interval setting, the complete indifference is represented by $[0,0]$, but we need to introduce here a more flexible representation of the neutral uncertainty. The reason behind this lies in interval arithmetic. Given two intervals $z_1,$ $z_2 \in \ci(\R)$, their difference should be all possible differences between elements in $z_1$ and $z_2$, i.e. 
\begin{equation*}
    z_1-z_2=\{x_1-x_2\tq x_1\in z_1,\quad x_2\in Z_2\}=[z_1^--z_2^+,z_1^+-z_2^-].
\end{equation*}

As a consequence, we have the following.
\begin{equation*}
    z-z=[z^--z^+,\,z^+-z^-]=[-\ell(z),\,\ell(z)],
\end{equation*}
which can be interpreted as an uncertain extension of the neutral value in a bipolar scale. However, we must emphasize that, while on a classical bipolar scale there is a single value representing neutrality, in the interval bipolar scale $\ci(\R)$ there are infinitely many interval values that are candidates to represent the idea of neutrality. Therefore, we consider the set.
\begin{equation*}
    \cn=\{[-\varepsilon,\,\varepsilon] \tq \varepsilon\geq0\}
\end{equation*}

Observe that $\cn$ is exactly the set of the self-differences between the elements in $\ci(\R)$. Later, we will further discuss which of these possible intervals is most suitable to represent the notion of neutrality. In addition, it is important to note that the opposite element for $z\in\ci(\R)$ is given by $[-z^+,\,-z^-]\in\ci(\R)$, that is, an element such that $z+(-z)=[-\ell(z),\,\ell(z)]$. For this reason, in the remainder of the manuscript, we extend the notation and define the opposed element as $-z=[-z^+,\,-z^-]$.

From the IPR point of view, we should keep in mind that comparing an object with itself should always produce a value within $\cn$, i.e., $z_{ii}\in\cn, \pt i=1,...,n$. Furthermore, from the previous discussion,  we can conclude that the reciprocity condition for IPRs may be stated as follows. 
\begin{equation*}
    z_{ij}+z_{ji}=[-\ell(z_{ij}),\,\ell(z_{ij})], \pt i,j=1,...,n,
\end{equation*}
or, equivalently, $z_{ij}=-z_{ji}$. Now, we can state the formal definition for IPR.
\begin{definition}[Interval Preference Relation (IPR)] Given an interval-valued matrix $Z\in\cmi$, we say that $Z$ is an IPR if it satisfies the reciprocity condition:
\begin{equation*}
    z_{ij}=-z_{ji}, \pt i,j=1,...,n.
\end{equation*}
\end{definition}

Note that, in particular, an IPR satisfies that $z_{ii}=-z_{ii}$ and thus $z_{ii}\in\cn\pt i=1,...,n$. Still, we need to clarify the notion of consistent IPR. In this case, we pursue the ideal case in which a single column or row allows for reconstructing the entire matrix by also generalizing the classical definitions for consistency to the interval settings. A first reasonable attempt could be defining consistency as $z_{ij}=z_{ik}+z_{kj}, \pt i,j,k=1,...,n$. However, this definition may fail even with simple examples. For instance, a particular case of this definition is $z_{11}=z_{11}+z_{11}$, which forces $z_{11}=[0,0]$, and thus we lose all the potential of providing a good generalization. Instead, we provide the following definition 

\begin{definition}[Consistency]
An IPR $Z$ is said to be consistent with respect to the neutral element $u\in\cn$ if it satisfies:
\begin{equation*}
    z_{ij}+u=z_{ik}+z_{kj}, \pt i,j,k=1,...,n.
\end{equation*}
\end{definition}
Note that we have inserted in the definition a neutral element $u\in\cn$ that indeed generalizes the role of $0$ to the interval case. Of course, if $u=[0,0]$ and $\ell(z_{ij})=0 \pt i,j=1,...,n$, we obtain a definition of consistency for the real number case. In the next section, we will show that such a definition is equivalent to the definitions of consistency given for other classical PRs. Let us discuss some properties that can be deduced directly from this definition.
\begin{proposition}
    Let $Z$ denote a consistent IPR with respect to $u\in\cn$. Then $z_{ii}=u\pt i=1,...,n$.
\end{proposition}
\begin{proof}
    Since $z_{ii}+u=z_{ii}+z_{ii}$, this implies that
    \begin{gather*}
        z_{ii}^-+u^-=z_{ii}^-+z_{ii}^-\\
        z_{ii}^++u^+=z_{ii}^++z_{ii}^+
    \end{gather*}
    from where $z_{ii}=u, \pt i=1,...,n$.
\end{proof}

From this result, we find that in a consistent IPR, all the interval elements in the main diagonal have the same length. Consequently,
\begin{corollary}\label{coro:length}
    If $Z$ denotes a consistent IPR with respect to $u\in\cn$, then 
    \begin{equation*}
        \ell(z_{ij})=\ell(u), \pt i,j=1,...,n
    \end{equation*}
\end{corollary}
\begin{proof}
    Since $z_{ii}+u=z_{ij}+z_{ji}, \pt i,j=1,...,n$, then $2u=[-\ell(z_{ij}),\,\ell(z_{ij})], \pt i,j=1,...,n$, from which $u^-=u^+=\ell(z_{ij}), \pt i,j=1,...,n$.
\end{proof}
In terms of the units, this means that the consistent scenario only admits interval values with the same units of uncertainty. Of course, in practice, it is hard to obtain this ideal scenario, in the same way that obtaining a perfectly consistent matrix is unlikely. This is a practical limitation that we address later through a socio-technical approach based on optimization models.

At this stage, let us prove the main result concerning consistent IPRs, which is indeed a representation theorem.
\begin{theorem}[Representation of consistent IPRs]\label{teor:repr}
    An IPR $Z$ is consistent with respect to $u\in\cn$ if and only if there exists an $n$-tuple of priority interval values $v_1,...,v_n\in\ci(\R)$ such that $z_{ij}+u=v_i-v_j, \pt i,j=1,...,n$.
\end{theorem}
\begin{proof}
    Let us assume first that $Z$ is consistent with respect to $u\in\cn$. Then, we can consider the $n$-tuple given by $v_k=z_{k,1}, \pt k=1,...,n$.
    Then, the consistency and reciprocity properties imply that
    \begin{equation*}
        z_{ij}+u=z_{i1}+z_{1j}=z_{i1}-z_{j1}=v_i-v_j, \pt i,j=1,...,n.
    \end{equation*}
    Conversely, let us assume that there exists an $n$-tuple of priority interval values $v_1,...,v_n\in\ci(\R)$ such that $z_{ij}+u=v_i-v_j, \pt i,j=1,...,n$. Then, using that $2u=z_{kk}+u=v_k-v_k, \pt k=1,...,n$, we have that
    \begin{equation*}
         z_{ik}+z_{kj}+2u=z_{ik}+u+z_{kj}+u=v_i-v_k+v_k-v_j=v_i-v_j+v_k-v_k=z_{ij}+u+z_{kk}+u=z_{ij}+3u
    \end{equation*}
    which implies that $z_{ik}+z_{kj}=z_{ij}+u$.
\end{proof}
The previous result guarantees that, as in the classical real-valued case, a consistent IPR contains the same information as one of its columns or rows.

\section{Interval-valued value scales from IPRs}\label{sec:scales}
 \noindent Now, let us emphasize here another consequence of the previous representation theorem. If we assume that the lowest possible performance value is given by comparing the worst object with itself, i.e., $v_n=z_{n,n}$, it is clear that the comparisons of the other objects with this one may provide a measure of the global performance of each object. Consequently, we can use the previous result to obtain an interval-valued value function from the preferences. 
\begin{proposition}
Suppose that $Z$ is a consistent IPR such that the compared objects are ordered from best to worst $1\succeq 2\succeq....\succeq{n}$, i.e., $z_{i(i+1)}\geq_0 u, \pt i=1,...,n-1$. Then, the interval $n$-tuple given by $v_k=z_{kn}, \pt k=1,...,n$ satisfies $v_k\geq_0 v_{k+1}, \pt i=1,...,n$.
\end{proposition}
\begin{proof}
For $v_k=z_{kn}, \pt k=1,...,n$, we have
\begin{equation*}
    v_k+u=z_{kn}+u= z_{k(k+1)}+z_{(k+1)n}=  z_{k(k+1)}+v_{k+1}.
\end{equation*}
Consequently, we have the following equivalences:
\begin{gather*}
    v_k+u=z_{k(k+1)}+v_{k+1} \iff \left\{\begin{array}{c}
         v_k^++u^+=z_{k(k+1)}^++v_{k+1}^+\\
          v_k^-+u^-=z_{k(k+1)}^-+v_{k+1}^-
    \end{array}\right.
    \iff \left\{\begin{array}{c}
         v_k^+-v_{k+1}^+=z_{k(k+1)}^+-u^+\\
          v_k^--v_{k+1}^-=z_{k(k+1)}^--u^-
    \end{array}\right.
\end{gather*}
This allows deducing that $v_k^+-v_{k+1}^+\geq 0$ and $v_k^--v_{k+1}^-\geq 0$, which is equivalent to the lattice order $\leq_0$.
\end{proof}
We remark here that we get the neutral element of the scale as the worst performance level because we have used as a reference the value $z_{n,n}$. However, any of the intervals in the $n$-th row of the matrix could have been chosen as a reference value for the worst performance level. By doing so, the performance value for the $i$-th object would be given by the $i$-th entry of the corresponding column. However, doing this does not affect the differences between different performance levels of the objects. 

Further, we emphasize that the values $v_1,...,v_n$ determine an interval-valued value scale. In this sense, we know that the monotonicity condition $v_1\geq_0 v_2\geq_0...\geq_0v_n$ is guaranteed (under the hypothesis of the previous result). This property is essential from the operational point of view since it implies that these quantities keep their order for every admissible order \cite{SANTANA202044, WOS:000317886300005}.

\begin{example}
    We highlight here that the definitions given before do not depend on a scale factor. Suppose that the original information is the number of units between consecutive objects that are ordered from best to worst. For example, $z_{12}=[4,6],z_{23}=[1,3]$, and $z_{34}=[2,4]$. Then, for the neutrality value $[-1,1]$, we obtain:
    \begin{equation*}
        \begin{pmatrix}
        [-1,1] & [4,6]   & [6,8]   & [9,11]  \\
        [-6,-4] & [-1,1]  & [1,3]   & [4,6]  \\
        [-8,-6]  & [-3,-1] & [-1,1]  & [2,4]  \\
        [-11,-9]  & [-6,-4]  & [-4,-2] & [-1,1] \\
        \end{pmatrix}
    \end{equation*} 
    And thus, the interval-valued value function is $v_1=[9,11]$, $v_2=[4,6]$, $v_3=[2,4]$ and $v_4=[-1,1]$. The detailed steps regarding the computation of these values is shown in Section \ref{sec:example}.
\end{example}
From this example, we can see that the conditions that determine consistent IPR are satisfied. However, it is natural to wonder about choosing a scale factor to normalize the matrix. From now on, and following a previous discussion, we consider that the lowest performance level will be $u$. This is only an assumption that we make here because $u$ generalizes the notion of $0$, and thus it is a good candidate for the sake of readability. The following result 11gives an idea of the relationship between different performance levels:
\begin{proposition}
    Suppose that $Z$ is a consistent IPR with respect to $u\in\cn$ such that the compared objects are ordered from best to worst $1\succeq 2\succeq....\succeq{n}$, i.e., $z_{i(i+1)}\geq_0 u, \pt i=1,...,n-1$. Then, the interval $n$-tuple given by $v_k=z_{kn}, \pt k=1,...,n$ satisfies
    \begin{equation*}
        v_k+(n-k-1)u=\sum_{i=k}^{n-1}z_{i(i+1)}, \pt k=1,...,n-1.
    \end{equation*}
\end{proposition}
\begin{proof}
    Let us proceed by reverse induction on $k=n-1,...,1$. The case $k=n-1$ is true since $v_{n-1}=z_{(n-1)n}$. Let us assume that the following inequalities hold for a fixed value $k\in\{2,...,n-1\}$:
    \begin{equation*}
        v_j+(n-j-1)u=\sum_{i=j}^{n-1}z_{i(i+1)}, \pt j=k,...,n-1
    \end{equation*}
    Then,     
    \begin{gather*}
       v_{k-1}+(n-(k-1)-1)u=z_{(k-1)n}+u+(n-k-1)u\\=z_{(k-1)k}+z_{kn}+(n-k-1)u=z_{(k-1)k}+\sum_{i=k}^{n-1}z_{i(i+1)}=\sum_{i=k-1}^{n-1}z_{i(i+1)}.
    \end{gather*}
\end{proof}

From the previous results, we can identify that the highest performance level is given by  $v_1=z_{1n}$, which satisfies 
\begin{equation*}
    \begin{cases}        
     {\displaystyle v_1^-=\sum_{i=1}^{n-1}z_{i(i+1)}^--(n-2)u^-}\\    
     {\displaystyle v_1^+=\sum_{i=1}^{n-1}z_{i(i+1)}^+-(n-2)u^+}\\    
    \end{cases}.
\end{equation*}
Notice that the quantities $C^-:=\sum_{i=1}^{n-1}z_{i(i+1)}^-$ and $C^+:=\sum_{i=1}^{n-1}z_{i(i+1)}^+$ represent the lowest and highest number of units used to make the pairwise comparisons among consecutive elements. These equations may be rewritten in terms of length $\ell(u)$, resulting in:
\begin{equation*}
    \begin{cases}        
    {\displaystyle v_1^-=C^-+(n-2)\ell(u)/2}\\    
    {\displaystyle v_1^+=C^+-(n-2)\ell(u)/2}\\    
    \end{cases}.
\end{equation*}

Therefore, the middle point of $v_1$ is exactly the average between the maximum and minimum number of units $C:=\frac{C^++C^-}{2}$. Consequently, if we use $\frac{1}{C}$ as a scale factor, we obtain that $\frac{1}{C}v_1$ (which corresponds to the maximum performance level) is an interval around $1$, i.e. $\frac{1}{C}v_1^-\leq 1\leq \frac{1}{C}v_1^+$. We use the process described here as the standard normalization method since it is in line with the classical vision. Of course, in our interval-valued version, the minimum performance level will be a generalized/interval version of $0$ ($\frac{1}{C}u$) and the highest performance level will be a generalized/interval version of $1$ ($\frac{1}{C}v_1$). Note that the scale factor $\frac{1}{C}>0$ is completely compatible with all the properties described in this section, assuming that, after normalizing, the role of the neutral element is played by $\frac{1}{C}u$ instead of $u$.

\section{IPR and classical preference relations}\label{sec:classical}
\noindent This section aims to show that our framework generalizes many of the classical preference relations in the literature. First, let us present an alternative version of our former Representation Theorem.
\begin{theorem}[Representation Theorem (revisited)]\label{teor:repr-crisp}
    Let $Z$ denote a consistent IPR with respect to $u\in\cn$. Then $Z$ is consistent if and only if there exists a consistent IPR $X$ with respect to $[0,0]$ such that $z_{ij}=\left[x_{ij}-\frac{\ell(u)}{2},x_{ij}+\frac{\ell(u)}{2}\right], \pt i,j=1,...,n$. Further, $v_1,...,v_n$ are interval-valued priority values for $Z$ if and only if $\nu_1=\frac{v_1^-+v_1^+}{2},...,\nu_n=\frac{v_n^-+v_n^+}{2}$ are priority values for $X$.
\end{theorem}
\begin{proof}
    Let $Z$ be a consistent IPR with respect to $u\in\cn$. Then, by Corollary \ref{coro:length}, we know that all the entries have the same length $\ell(u)$. Let us define $x_{ij}=\frac{z_{ij}^++z_{ij}^-}{2}\pt i,j=1,...,$ and consider the IPR $X$ whose items are of the shape $[x_{ij},x_{ij}]$, which satisfies $\ell(x_{ij})=0\pt i,j=1,...,n$ and also $z_{ij}=\left[x_{ij}-\frac{\ell(u)}{2},x_{ij}+\frac{\ell(u)}{2}\right], \pt i,j=1,...,n$. Further, the consistency of $Z$ implies
    \begin{gather*}
        z_{ij}^++u^+=z_{ik}^++z_{kj}^+\\
        z_{ij}^+-u^+=z_{ik}^-+z_{kj}^-
    \end{gather*}
    which can be rewritten in terms of $X$ as 
    \begin{gather*}
        x_{ij}+\frac{\ell(u)}{2}+u^+=x_{ik}+\frac{\ell(u)}{2}+x_{kj}+\frac{\ell(u)}{2}\\
        x_{ij}-\frac{\ell(u)}{2}-u^-=x_{ik}-\frac{\ell(u)}{2}+x_{kj}-\frac{\ell(u)}{2}
    \end{gather*}
    or equivalently $ x_{ij}=x_{ik}+x_{kj}$ which is the consistency of $X$ with respect to $[0,0]$. From the previous reasoning, we can also conclude the reciprocal statement, i.e., if if $X$ is a consistent IPR with respect to $[0,0]$, the IPR $Z$ constructed by $z_{ij}=x_{ij}+\left[-\frac{\ell(u)}{2},\frac{\ell(u)}{2}\right]$ is also consistent. A direct application of Theorem \ref{teor:repr} allows us to obtain the last statement.
\end{proof}
We emphasize here that an IPR with respect to $[0,0]$ is nothing but a real-valued matrix $X$ that satisfies $x_{ij}=x_{ik}+x_{kj}, \pt i,j,k=1,...,n$. In particular, this definition implies $x_{ii}=0, \pt i=1,...,n$ and $x_{ij}=-x{ji}, \pt i,j=1,...,n$. Notice that this is exactly the notion of pairwise consistent table utilized in \cite{corrente2021} that was integrated within the Deck of Cards framework.

Also, let us recall that consistent FPRs \cite{Orlovski1978} are characterized by the relation $y_{ij}+\frac{1}{2}=y_{ik}+y_{kj}, \pt i,j,k=1,...,n$. Note that the transformation $f(y)=y-\frac{1}{2}$ allows remapping the FPR $Y$ onto an interval IPR. Further, it can be easily shown that the FPR $Y$ is consistent if and only if $f(Y)$ is consistent with respect to $[0,0]$. 

We could also consider the MPRs introduced by Saaty \cite{Saaty2008}, whose consistency condition was given by $a_{ij}=a_{ik}a_{kj}, \pt i,j,k=1,...,n$. In this case, we can consider the transformation $f(a)=\log_9(a)$ to remap MPRs into IPRs. Further, it is also possible to prove that an MPR $A$ is consistent if and only if the IPR $f(A)$ is consistent with respect to $[0,0]$.

Notice that, in this sense, the framework here proposed generalizes all these existing approaches by also offering a common ground on which all these notions are equivalent. Additionally, the procedure proposed here to obtain value scales can also be applied to compute value scales for all the previous types of PRs.

Furthermore, our framework also provides a methodology to work with interval-valued versions of FPRs, MPRs, and DoC. To do so, we only need to apply the above-described transformation to embed the corresponding preference relation into the IPR here defined, compute the interval-valued value scales, and apply the inverse transformation to the result. We further explore this process in Section \ref{sec:example} for the case of the DoC.

\section{Optimization models for practical scenarios}\label{sec:opt}
\noindent Until now, we have described how to work with IPRs which are consistent with respect to a given interval $u\cn$. However, we acknowledge that in practice, eliciting a consistent IPR from the decision-maker may be hard in practice. For this reason, we devote this section to presenting several optimization models that may help the analyst to guide the decision-maker towards a consistent IPR. Let us start with the case in which the decision-maker provides a complete interval-valued pairwise comparison table $Z$ that may not be consistent. Our purpose is to identify for which value $u\in\cn$ we can find an IPR $\overline{Z}$ consistent with respect to $u$ and that is similar to the original matrix. Keeping this in mind, we propose the following optimization model:

\begin{equation}\label{model:l2}\tag{$M_{l_2}$}
    \begin{array}{rl}
        {\displaystyle \min_{\overline{z}_{ij}\in\ci(\R), u\in\cn}} &  {\displaystyle\frac{1}{2n^2}\sum_{i=1}^n\sum_{j=1}^n(z_{ij}^+-\overline{z}_{ij}^+)^2+(z_{ij}^--\overline{z}_{ij}^-)^2} \\
        \mbox{subject to:} & {\displaystyle \overline{z}_{ij}+u=\overline{z}_{ik}+\overline{z}_{kj}}, \; \pt i,j,k=1,...,n\\
    \end{array}
\end{equation}

Note that the objective function is based on the $l_2$ norm. However, we could use instead the $l_1$ or $l_\infty$ norms and obtain a linear programming model. In this paper, we keep the $l_2$ version because of the convenient analytical properties provided by the following result.

\begin{theorem}[Existence and uniqueness conditions for Model (\ref{model:l2})]\label{teor:opt}
    Let us assume that $Z$ is an IPR. Then, the optimization Model (\ref{model:l2}) is equivalent to:
    \begin{equation*}
    \begin{array}{rl}
        {\displaystyle \min_{\nu_1,...,\nu_n\in\R, \alpha\geq0}} &  {\displaystyle\frac{1}{2n^2}\sum_{i=1}^n\sum_{j=1}^n(z_{ij}^+-(\nu_i-\nu_j+\alpha))^2+(z_{ij}^--(\nu_i-\nu_j-\alpha))^2} \\
    \end{array}
\end{equation*}
Consequently, an optimal solution for \ref{model:l2} is given by
\begin{equation*}
    \begin{cases}
        {\displaystyle \nu_k^*=\frac{1}{n}\sum_{j=1}^n\nu_j+\frac{1}{n}\sum_{j=1}^nc_{kj}\pt k=1,...,n}\\
        {\displaystyle \alpha^*=\frac{1}{2n^2}\sum_{i=1}^n\sum_{j=1}^nl_{ij}}
    \end{cases}.
\end{equation*}
where $c_{ij}=\frac{z_{ij}^++z_{ij}^-}{2}$ and $l_{ij}=\ell(z_{ij})$. Furthermore, for a fixed $\mu>0$, the model 
   \begin{equation*}
    \begin{array}{rl}
        {\displaystyle \min_{\nu_1,...,\nu_n\in\R, \alpha\geq0}} &  {\displaystyle\frac{1}{2n^2}\sum_{i=1}^n\sum_{j=1}^n(z_{ij}^+-(\nu_i-\nu_j+\alpha))^2+(z_{ij}^--(\nu_i-\nu_j-\alpha))^2} \\
         \mbox{subject to:} & {\displaystyle \frac{1}{n}\sum_{i=1}^nv_i=\mu}\\
    \end{array}
\end{equation*}
has a unique solution
\begin{equation*}
    \begin{cases}
        {\displaystyle \nu_k^*=\mu+\frac{1}{n}\sum_{j=1}^nc_{kj}, \pt k=1,...,n}\\
        {\displaystyle \alpha^*=\frac{1}{2n^2}\sum_{i=1}^n\sum_{j=1}^nl_{ij}}
    \end{cases}.
\end{equation*}
which is also a solution for \ref{model:l2}.
\end{theorem}
\begin{proof}
    Using Theorem \ref{teor:repr-crisp}, since $\overline{Z}$ is a consistent IPR with respect to $u$, we can find $\nu_1,...,\nu_n$ such that $\overline{z}_{ij}+u=\left[\nu_{i}-\frac{\ell(u)}{2},\nu_i+\frac{\ell(u)}{2}\right]+\left[-\nu_j-\frac{\ell(u)}{2},-\nu_j+\frac{\ell(u)}{2}\right]$. Equivalently, 
    \begin{equation*}
        \begin{cases}
            \overline{z}_{ij}^++\frac{\ell(u)}{2}=\nu_i+\frac{\ell(u)}{2}+-\nu_j+\frac{\ell(u)}{2}\\
        \overline{z}_{ij}^--\frac{\ell(u)}{2}=\nu_{i}-\frac{\ell(u)}{2}-\nu_j-\frac{\ell(u)}{2}
        \end{cases}\iff \begin{cases}
            \overline{z}_{ij}^+=\nu_i-\nu_j+\frac{\ell(u)}{2}\\
        \overline{z}_{ij}^-=\nu_{i}-\nu_j-\frac{\ell(u)}{2}
        \end{cases}
    \end{equation*}
    and renaming $\alpha=\frac{\ell(u)}{2}$ gives us the equivalence between \ref{model:l2} and the first model in the statement of the theorem. Now, to compute an optimal analytical solution for such a model, let us consider $f_{ij}(\nu,\alpha)=(z_{ij}^+-(\nu_i-\nu_j+\alpha))^2+(z_{ij}^--(\nu_i-\nu_j-\alpha))^2$ and compute the partial derivatives. Notice that, for $k\in\{1,...,n\}$ we have the following cases:
    \begin{itemize}[label={--}]
        \item $k=i$, $k=j$. $\partial_{\nu_k}f_{ij}(\nu,\alpha)=0$
        \item $k=i$, $k\neq j$. $\partial_{\nu_k}f_{kj}(\nu,\alpha)=-2(z_{kj}^+-(\nu_k-\nu_j+\alpha))-2(z_{kj}^--(\nu_k-\nu_j-\alpha))=-2(z_{kj}^++z_{kj}^--2(\nu_k-\nu_j))$.
        \item $k\neq i$, $k=j$. $\partial_{\nu_k}f_{ik}(\nu,\alpha)=2(z_{ik}^+-(\nu_i-\nu_k+\alpha))+2(z_{ik}^--(\nu_i-\nu_k-\alpha))=2(z_{ik}^++z_{ik}^--2(\nu_i-v_k))$.
        \item $k\neq i$, $k\neq j$. $\partial_{\nu_k}f_{ij}(\nu,\alpha)=0$.
    \end{itemize}
    Further, $\partial_{\alpha}f_{ij}(\nu,\alpha)=-2(z_{ij}^+-(\nu_i-\nu_j+\alpha))+2(z_{ij}^--(\nu_i-\nu_j-\alpha))=2(z_{ij}^--z_{ij}^++2\alpha)$. Now. we can derive the partial derivatives of the objective funtion. If we define as follows
    \begin{equation*}
        F(\nu,\alpha)=\frac{1}{2n^2}\sum_{i=1}^n\sum_{j=1}^n(z_{ij}^+-(\nu_i-\nu_j+\alpha))^2+(z_{ij}^--(\nu_i-\nu_j-\alpha))^2=\frac{1}{2n^2}\sum_{i=1}^n\sum_{j=1}^nf_{ij}(\nu,\alpha),
    \end{equation*}
    we obtain
    \begin{gather*}        \partial_{\nu_k}F(\nu,\alpha)=\partial_{\nu_k}\frac{1}{2n^2}(\sum_{i=1}^n\sum_{j=1}^nf_{ij}(\nu,\alpha))=\partial_{\nu_k}\frac{1}{2n^2}(\sum_{i=k}\sum_{j=1}^nf_{ij}(\nu,\alpha)+\sum_{i\neq k}\sum_{j=1}^nf_{ij}(\nu,\alpha))\\
        =\partial_{\nu_k}\frac{1}{2n^2}(\sum_{i=k}(\sum_{j=k}f_{ij}(\nu,\alpha)+\sum_{j\neq k}f_{ij}(\nu,\alpha))+\sum_{i\neq k}(\sum_{j=k}f_{ij}(\nu,\alpha)+\sum_{j\neq k}f_{ij}(\nu,\alpha)))\\
        =\frac{1}{2n^2}(\partial_{\nu_k}f_{kk}(\nu,\alpha)+\sum_{j\neq k}\partial_{\nu_k}f_{kj}(\nu,\alpha)+\sum_{i\neq k}\partial_{\nu_k}f_{ik}(\nu,\alpha)+\sum_{i\neq k}\sum_{j\neq k}\partial_{\nu_k}f_{ij}(\nu,\alpha))\\
        =\frac{1}{2n^2}(\sum_{j\neq k}\partial_{\nu_k}f_{kj}(\nu,\alpha)+\sum_{i\neq k}\partial_{\nu_k}f_{ik}(\nu,\alpha))=\frac{1}{2n^2}(\sum_{j\neq k}\partial_{\nu_k}f_{kj}(\nu,\alpha)+\partial_{\nu_k}f_{jk}(\nu,\alpha))\\
        =\frac{1}{2n^2}\sum_{j\neq k}(-2(z_{kj}^++z_{kj}^--2(\nu_k-\nu_j))+2(z_{jk}^++z_{jk}^--2(\nu_j-v_k))\\
        =\frac{1}{n^2}\sum_{j\neq k}(4(\nu_k-\nu_j)+z_{jk}^+-z_{kj}^++z_{jk}^--z_{kj}^-)=\frac{1}{n^2}\sum_{j\neq k}(4(\nu_k-\nu_j)+z_{jk}^++z_{jk}^-+z_{jk}^-+z_{jk}^+)\\
        =\frac{2}{n^2}\sum_{j\neq k}(2(\nu_k-\nu_j)+z_{jk}^++z_{jk}^-),
    \end{gather*}
    and also
    \begin{gather*}
        \partial_{\alpha}F(\nu,\alpha)=\partial_{\alpha}\frac{1}{2n^2}\sum_{i=1}^n\sum_{j=1}^nf_{ij}(\nu,\alpha)=\frac{1}{n^2}\sum_{i=1}^n\sum_{j=1}^n(z_{ij}^--z_{ij}^++2\alpha).
    \end{gather*}
Therefore,
\begin{equation*}
    \begin{cases}
        \partial_{\nu_k}F(\nu,\alpha)=0, \pt k=1,...,n\\
        \partial_{\alpha}F(\nu,\alpha)=0
    \end{cases}\iff
    \begin{cases}
        {\displaystyle \sum_{j=1}^n(\nu_k-\nu_j+c_{jk})=0\pt k=1,...,n} \\
        {\displaystyle \sum_{i=1}^n\sum_{j=1}^n(2\alpha-l_{ij})=0}
    \end{cases}
\end{equation*}
Equivalently, 
\begin{equation*}
    \begin{cases}
        {\displaystyle \nu_k=\frac{1}{n}\sum_{j=1}^n\nu_j+\frac{1}{n}\sum_{j=1}^nc_{kj}, \pt k=1,...,n}\\
        {\displaystyle \alpha=\frac{1}{2n^2}\sum_{i=1}^n\sum_{j=1}^nl_{ij}}
    \end{cases}.
\end{equation*}

Finally, let us focus on the last optimization model. Let us proceed via Lagrange multipliers to find an optimal solution for it. Therefore, let us consider $L(\nu,\alpha,\lambda)=F(\nu,\alpha)+\lambda(\mu-\frac{1}{n}\sum_{i=1}^nv_i)$. Consequently,
\begin{gather*}
     \partial_{\nu_k}L(\nu,\alpha,\lambda)= \partial_{\nu_k}F(\nu,\alpha) -\lambda\frac{1}{n}v_k, \pt k=1,...,n\\
     \partial_{\alpha}L(\nu,\alpha,\lambda)=\partial_{\alpha}F(\nu,\alpha)\\
    \partial_{\lambda}L(\nu,\alpha,\lambda)=\mu-\frac{1}{n}\sum_{i=1}^nv_i
\end{gather*}
and the partial derivatives are equal to zero if and only if 
\begin{equation*}
    \begin{cases}
        {\displaystyle \frac{4}{n}\sum_{j=1}^n((\nu_k-\nu_j)+c_{jk})= \lambda v_k, \pt k=1,...,n}\\
        {\displaystyle \sum_{i=1}^n\sum_{j=1}^n(2\alpha-l_{ij})=0} \\
        {\displaystyle \mu=\frac{1}{n}\sum_{i=1}^nv_i}
    \end{cases}
\end{equation*}
Notice that
\begin{gather*}
    \lambda\mu={\lambda}\sum_{k=1}^nv_k=\sum_{k=1}^n \frac{4}{n}\sum_{j=1}^n((\nu_k-\nu_j)+c_{jk})=\frac{4}{n}\sum_{k=1}^nn\nu_k \sum_{j=1}^n(c_{jk}-\nu_j)=\frac{4}{n}\sum_{k=1}^n(n\nu_k+ \sum_{j=1}^nc_{jk}-\mu)\\=4(n\mu +\sum_{k=1}^n\sum_{j=1}^nc_{jk}-n\mu)=0.
\end{gather*}
From here, $\lambda=0$ and consequently
\begin{equation*}
    \begin{cases}
        {\displaystyle \nu_k^*=\mu+\frac{1}{n}\sum_{j=1}^nc_{kj}, \pt k=1,...,n}\\
        {\displaystyle \alpha^*=\frac{1}{2n^2}\sum_{i=1}^n\sum_{j=1}^nl_{ij}}
    \end{cases}.
\end{equation*}
is an optimal solution for the second optimization model that is also optimal for the first one. Further, given two solutions of this second optimization model $(\nu^1,\alpha^1)$ and $(\nu^2,\alpha^2)$, since they must satisfy the previous equations, it must be
\begin{equation*}
    \begin{cases}
        {\displaystyle \nu_k^1=\mu+\frac{1}{n}\sum_{j=1}^nc_{kj}=\nu_k^2, \pt k=1,...,n}\\
        {\displaystyle \alpha^1=\frac{1}{2n^2}\sum_{i=1}^n\sum_{j=1}^nl_{ij}=\alpha^2}
    \end{cases}.
\end{equation*}
and thus the solution for that model is unique.
\end{proof}

Let us hihglight here that after solving the Model (\ref{model:l2}), we obtain both the neutral interval $u=[-\alpha,\alpha]$ and interval-valued value scale $[v_1-\alpha,v_1+\alpha],...,[v_n-\alpha,v_n+\alpha]$ that allow generating the IPR $\overline{z}_{ij}=[v_i-v_j-\alpha,v_i-v_j+-\alpha]\pt i,j=1,...,n$ which is consistent with respect to $u$ and which is closest to the original preference table $Z$. In practice, this information will be shared with the decision-maker in order to improve the consistency of his/her opinions.

Further, keep in mind that if the preference table is given in terms of units, it is sensitive to expect that $\ell(u)$ is a non-negative integer. At this stage, one possible strategy consists of including such a restriction as a constraint on the previous Model (\ref{model:l2}) (or any of the variants $l_1$ or $l_\infty$), obtaining an integer programming problem. However, here we provide a more user-friendly approach that is more in line with the previous result, using that the value $\ell(u)\in\Z$ could be approximately inferred from the context. For a fixed $u=[-\alpha,\alpha]$ ($\alpha>0$), let consider the following optimization problem:

\begin{equation}\label{model:l2-ufix}\tag{$M_{l_2},u$}
    \begin{array}{rl}
        {\displaystyle \min_{\overline{z}_{ij}\in\ci(\R)}} &  {\displaystyle\frac{1}{2n^2}\sum_{i=1}^n\sum_{j=1}^n(z_{ij}^+-\overline{z}_{ij}^+)^2+(z_{ij}^--\overline{z}_{ij}^-)^2} \\
        \mbox{subject to:} & {\displaystyle \overline{z}_{ij}+u=\overline{z}_{ik}+\overline{z}_{kj}}, \; \pt i,j,k=1,...,n\\
    \end{array}
\end{equation}
which is exactly the same as the former Model (\ref{model:l2}) but removing $u$ from the set of variables. Using a reasoning analogous to Theorem \ref{teor:opt}, we can obtain that 
\begin{equation*}
    \nu_k^*=\frac{1}{n}\sum_{j=1}^n\nu_j+\frac{1}{n}\sum_{j=1}^nc_{kj}, \pt k=1,...,n
\end{equation*}
determins a solution for Model (\ref{model:l2-ufix}) by taking $\overline{z}_{ij}=[v_i-v_j-\alpha,\; v_i-v_j+\alpha],\, i,j=1,...,n$. Consequently, instead of solving an integer programming problem, we may solve the previous model for some values of $\alpha$ and show the results to the decision-maker in order to guide him/her towards a consistent preference table.

\section{Illustrative Examples}\label{sec:example}
\noindent In this section, we present two examples regarding the computation of the interval-valued value scale. Both examples are based on the Deck of Cards Method \cite{corrente2021} and they cover both the case in which the provided intervals have the same length and the case with different lengths.
\subsection{Consistent pairwise table}

\begin{enumerate}
    \item  Suppose that we have the following four objects, totally ordered from the best to the worst: 
    \[ 
        \ell_1 \succeq \ell_2 \succeq \ell_3 \succeq \ell_4. 
    \]
    \item Insert an uncertain number of blank cards in between them. The blank cards represent the more os less attractiveness between two objects.
    \[
        \ell_1 \; [3,5] \; \ell_2 \; [0,2] \; \ell_3 \; [1,3] \; \ell_4. 
    \]
    \item The following incomplete pairwise table can be as follows:
    \begin{equation*}
        \begin{pmatrix}
          *   & [3,5] &       &        \\
             &   *    & [0,2] &        \\
             &       &    *   & [1,3]  \\
             &       &       &      *  \\
        \end{pmatrix}
    \end{equation*} 
    \item Zero blank cards does not mean the objects have the same value, but that this value is minimum, i.e., equal to the unit. Thus, the unit pairwise is built by summing one to each value of the previous table:
    \begin{equation*}
        \begin{pmatrix}
          *   & [4,6] &       &        \\
             &  *     & [1,3] &        \\
             &       &     *  & [2,4]  \\
             &       &       &     *   \\
        \end{pmatrix}
    \end{equation*} 
    \item By applying the consistency condition with respect to $u=[-1,1]$, we obtain the table:
    \begin{equation*}
        \begin{pmatrix}
            [-1,1] & [4,6] & [6,8] & [9,11] \\
             &  [-1,1]     & [1,3] & [4,6]  \\
             &       &  [-1,1]     & [2,4]  \\
             &       &       &   [-1,1]     \\
        \end{pmatrix}
    \end{equation*} 
    \item Now, we obtain a normalized value scale by dividing the last column by the average number of units, i.e., 
    \begin{equation*}
        \frac{(4+1+2)+(6+3+4)}{2}=10
    \end{equation*}
    which yields
    \begin{gather*}
        v_1\approx[0.9,1.1]\\
        v_2\approx[0.4,0.6]\\
        v_3\approx[0.2,0.4]\\
        v_4\approx[-0.1,0.1]\\
    \end{gather*}
    Note that these values satisfy the monotonicity condition $v_1\geq_0v_2\geq_0v_3\geq_0v_4$.
\end{enumerate}

\subsection{Inconsistent pairwise table}
\noindent In this case, we follow the guidelines of Section \ref{sec:opt} and provide a methodology to convert the original information given by the decision-maker into intervals of the same length before applying the computational rules to obtain the interval-valued value scale. To do so, let us consider the following optimization model, which provides intervals of the same length that are the closest to the intervals provided by the decision-maker:

\begin{equation}\label{model:DoCl2}\tag{$M_{DoC-l_2}$}
    \begin{array}{rl}
        {\displaystyle \min_{\overline{z}_{i(i+1)}\in\ci(\R), u\in\cn}} &  {\displaystyle\frac{1}{2(n-1)}\sum_{i=1}^{n-1}(z_{i(i+1)}^+-\overline{z}_{i(i+1)}^+)^2+(z_{i(i+1)}^--\overline{z}_{i(i+1)}^-)^2} \\
        \mbox{subject to:} & {\displaystyle \overline{z}_{i(i+1)}^+-\overline{z}_{i(i+1)}^-=\ell(u)}; \; \pt i=1,...,n-1\\
        &{\displaystyle \overline{z}_{i(i+1)}\geq_0 u}; \pt i=1,...,n-1\\
    \end{array}
\end{equation}
Similarly to the discussion in Section \ref{sec:opt}, we use the convention $u=[-\alpha,\alpha]$, to obtain the equivalent model
\begin{equation*}
    \begin{array}{rl}
        {\displaystyle \min_{\overline{z}_{i(i+1)}^+,\overline{z}_{i(i+1)}^-\in\R, \alpha\geq0}} &  {\displaystyle\frac{1}{2(n-1)}\sum_{i=1}^{n-1}(z_{i(i+1)}^+-\overline{z}_{i(i+1)}^+)^2+(z_{i(i+1)}^--\overline{z}_{i(i+1)}^-)^2} \\
        \mbox{subject to:} & {\displaystyle \overline{z}_{i(i+1)}^+-\overline{z}_{i(i+1)}^-=2\alpha}; \; \pt i=1,...,n-1\\
        &{\displaystyle \overline{z}_{i(i+1)}^+\geq\alpha}; \pt i=1,...,n-1\\
        &\alpha\geq0
    \end{array}
\end{equation*}
whose optimal solution is given by 
\begin{gather*}
    \alpha^*=\frac{1}{2(n-1)}\sum_{i=1}^{n-1}\ell(z_{i(i+1)})\\    (\overline{z}^{+}_{i(i+1)})^*=\frac{z_{i(i+1)}^++z_{i(i+1)}^-}{2}+\alpha^*\\    (\overline{z}^{-}_{i(i+1)})^*=\frac{z_{i(i+1)}^++z_{i(i+1)}^-}{2}-\alpha^*
\end{gather*}

\begin{enumerate}
    \item  As in the previous case, the starting point is four objects, totally ordered from the best to the worst: 
    \[ 
        \ell_1 \succeq \ell_2 \succeq \ell_3 \succeq \ell_4. 
    \]
    \item Now,  an uncertain number of blank cards is inserted between them, representing the more or less attractiveness between consecutive objects.
    \[
        \ell_1 \; [3,5] \; \ell_2 \; [0,2] \; \ell_3 \; [1,4] \; \ell_4. 
    \]
    \item Now, the pairwise comparison table is built by summing one to each value of the previous table:
    \begin{equation*}
        \begin{pmatrix}
          *   & [4,6] &       &        \\
             &    *   & [1,3] &        \\
             &       &   *    & [2,5]  \\
             &       &       &     *   \\
        \end{pmatrix}
    \end{equation*} 
    \item Note that in this case, the length of the given intervals is not the same. Before we proceed, the analyst must discuss with the decision-maker in order to obtain consistent information. To guide the analyst in this case, we may use the optimization Model (\ref{model:DoCl2}) introduced earlier. In this case, we obtain that the optimal values are:
    \begin{gather*}
        \alpha^*=\frac{2+2+3}{2\cdot3}=\frac{7}{6}\approx 1.166\\
        \overline{z}_{12}=\left[\frac{4+6}{2}-\alpha^*,\frac{4+6}{2}+\alpha^*\right]=[5-\alpha^*,5+\alpha^*]\approx[3.833,6.166]\\
        \overline{z}_{23}=\left[\frac{1+3}{2}-\alpha^*,\frac{1+3}{2}+\alpha^*\right]=[2-\alpha^*,2+\alpha^*]\approx[0.833,3.166]\\
        \overline{z}_{34}=\left[\frac{2+5}{2}-\alpha^*,\frac{2+5}{2}+\alpha^*\right]=[3.5-\alpha^*,3.5+\alpha^*]\approx[2.333,4.666]
    \end{gather*}
    At this stage, the analyst must double-check with the decision-maker if he/shefeels that these values represent his/her preferences. If not, further modifications should be made with the cards until an agreement is reached. Let us assume here that he/she feels that these values do represent his/her opinions. In that case, we may proceed as in the previous case to derive the value scale. 
    
    \item By applying the conistency condition with respect to $u=[-1.166,1.166]$, we obtain the table:
    \begin{equation*}
        \begin{pmatrix}
             [-1.166,1.166]& [3.833,6.166] & [5.333,8.166] & [9.333,11.666] \\
             &  [-1.166,1.166]     & [0.8333,3.166] & [4.333,6.666]  \\
             &       &  [-1.166,1.166]     & [2.333,4.666]  \\
             &       &       & [-1.166,1.166]       \\
        \end{pmatrix}
    \end{equation*} 
    \item Now, we may derive a normalized value scale by dividing the last column by the average number of units, i.e., 
    \begin{equation*}
        \frac{(4+1+2)+(6+3+5)}{2}=10.5
    \end{equation*}
    which yields
    \begin{gather*}
        v_1\approx[0.889,1.111]\\
        v_2\approx[0.421,0.634]\\
        v_3\approx[0.222,0.444]\\
        v_4\approx[-0.111,0.111]\\
    \end{gather*}
    Note that these values also satisfy the monotonicity condition $v_1\geq_0v_2\geq_0v_3\geq_0v_4$.
\end{enumerate}

\section{Conclusion}\label{sec:conclusion}
\noindent In this paper, we have presented a methodology to derive interval-valued value functions from interval-valued judgements. To do so, we have presented the notion of IPR and extended the classical notion of consistency to interval preference information. Furthermore, we have shown that if the interval judgments satisfy the proposed definition, interval-valued value functions can be univocally determined.

In addition, we show the unifying nature of our approach: by working with interval judgments, we recover as special cases several well-known models of preference representation. In particular, our framework generalizes additive transitivity in fuzzy preference relations, multiplicative consistency in Saaty’s preference relations, and the consistency definition in the Deck of Cards Method. This not only consolidates different paradigms under a common foundation, but also clarifies the role of consistency as the central principle linking them.

From a methodological perspective, we have proposed and solved several optimization models to verify whether a set of interval judgments satisfies consistency and to assist decision-makers in resolving inconsistencies when they appear. These models act both as diagnostic and corrective tools, making the framework suitable for interactive decision support.

The proposed theory and methods contribute to the development of well-defined and interpretable interval-valued scales, capable of capturing uncertainty without losing information in the process. Future research will focus on extending this framework to solve other classical decision-making problems, such as the construction of non-monotonic value functions or the proposal of other interval-arithmetic-based decision models, and exploring its application to real-world decision contexts where imprecision is unavoidable.

\section*{Acknowledgments}
\noindent The work is financed by Portuguese funds through the FCT – Foundation for Science and Technology under project UID/00097/2025 (CEGIST).

\vfill\newpage

\begingroup
\setstretch{1.0}  
\bibliographystyle{model2-names}
\bibliography{bibliography}
\endgroup
\end{document}